%
%
%
%


\magnification=\magstep1
\input amstex
\define\Bers{1}
\define\Hamen{2}
\define\Kgeo{3}
\define\Koranyi{4}
\define\Kcr{5}
\define\Kor{6}
\define\Lempert{7}
\define\Leme{8}
\define\Li{9}
\define\Pinch{10}
\define\Strich{11}
\define\Tang{12}   
\define\Teich{13}
\define\Vais{14}
\documentstyle{amsppt}
\topmatter
\nopagenumbers
\hyphenation {equi-valent}
\hyphenation {mani-fold}
\hyphenation {mini-mal}
\vsize=8.635truein
\hsize=6.55truein
\voffset0.3truein
\title
Regularity and extremality of quasiconformal homeomorphisms on cr 3-manifolds
\endtitle
\rightheadtext{quasiconformal homeomorphisms on cr manifolds}
\author
PUQI TANG
\endauthor
\abstract
This paper first studies the regularity of conformal homeomorphisms
on smooth locally embeddable strongly pseudoconvex CR manifolds. Then moduli
of curve families are used to estimate the maximal dilatations of
quasiconformal homeomorphisms. On certain CR 3-manifolds, namely, CR circle
bundles over flat tori, extremal quasiconformal homeomorphisms in some
homotopy classes are constructed. These extremal mappings have similar
behaviors to Teichm\"uller mappings on Riemann surfaces.  
\endabstract
\address
Department of Mathematics, Purdue University, West Lafayette, IN 47907, USA
\endaddress
\email
tang\@math.purdue.edu \newline \newline
\centerline{\boxed{\text{ Version 12.04.1993. Run at euclid.math.purdue.edu }}}
\endemail
\endtopmatter

\document

\heading
1. Introduction
\endheading

A contact manifold $M$ is a manifold of odd dimension with a non-integrable
distribution $HM$ of tangent hyperplanes. A Cauchy-Riemann (CR) manifold is
a contact manifold $M$ endowed a complex structure on the contact bundle $HM$.
Two CR manifolds are equivalent if there is a homeomorphism between them
which preserves both contact and CR structures. Generally, between any two
CR structures assigned on the same contact manifold, there may be no such
so-called CR homeomorphism between them. Therefore we consider those
homeomorphisms between CR manifolds which preserve the underlying contact
structures and distort the CR structures boundedly. They are called
quasiconformal homeomorphisms. In a class of homoemorphisms between two CR
manifolds, an extremal mapping is a quasiconformal homeomorphisms which
distorts the CR structures in a minimal way. This paper studies regularity
of quasiconformal homeomorphisms and extremal quasiconformal homeomorphisms
on smooth strongly pseudoconvex CR manifolds.

The notion of quasiconformal homeomorphisms is a new tool to study CR
structures as initiated by Kor\'anyi and Reimann. In this paper, we use
an analytic definition of quasiconformal homeomorphisms given in \cite{\Tang}
which is a generalization of the one given by Kor\'anyi and Reimann in
\cite{\Kcr}. We restrict ourselves to the 3-dimensional case here not only
because the notion of quasiconformality is not invariant under CR
transformations in higher dimensional cases, but also because 3-dimensional
CR structures are among the most interesting objects in the theory of
CR manifolds. We refer to \cite{\Lempert} for details about the second point.

Kor\'anyi and Reimann proved that $C^4$ conformal homeomorphisms on
Heisenberg groups must be smooth and CR (Theorem 8, \cite{\Kor}).
By applying a regularity theorem of weak CR mappings of Pinchuk and Tsyganov,
we generalize Kor\'anyi and Reimann's result to that a conformal
homeomorphism $f$ between two smooth, strongly pseudoconvex, locally embeddable
CR manifolds must be smooth and CR, if $f$ has $L^1_{loc}$ horizontal
derivatives (Theorem 2.3). Hence between such CR manifolds, quasiconformal
homeomorphisms with this weak regularity are actually ``quasi-CR".

To study the extremality of quasiconformal homeomorphisms is a global problem.
But our analytic definition of quasiconformality is proposed infinitesimally.
Therefore we need some global notion to describe the quasiconformality.
The one best fitting our later developments is the notion of moduli of curve
families. We prove that a $C^2$ diffeomorphism is quasiconformal if and only
if it preserves moduli of certain curve families up to a fixed bounded multiple
(Theorem 3.3). On the other hand, a homeomorphism satisfying this property
is absolutely continuous on lines (ACL) (Theorem 3.4).

Between CR circle bundles over flat tori, we construct extremal quasiconformal
homeomorphisms in certain homotopy classes (Theorem 4.2). There are two
transversal Legendrian foliations such that the extremal homeomorphism
constructed preserves these two foliations. More precisely, it is a stretching
by a constant factor along leaves of one foliation and a compressing by the
same factor along leaves of another foliation. This behavior is analogous to
those of Teichm\"uller mappings on Riemann surfaces (see \cite{\Teich} or
\cite{\Bers}). The generator $T$ of the circle action is transversal to the
contact bundle. In this transversal direction, the extremal mappings are
equivariant under the circle action.

But on an arbitrary CR 3-manifold with a transversal free circle action, an
extremal quasiconformal homeomorphism is not necessarily equivariant under
the circle action. Such CR manifolds are constructed in \cite{\Tang} so that
no extremal quasiconformal mapping between them is equivariant.

This work is heavily influenced by the theory of quasiconformal
homeomorphisms on Riemann surfaces and Teichm\"uller theory. Numerous
proofs in this paper are motivated by the proofs of the analogous facts
on Riemann surfaces. For example, the construction of the extremal
homeomorphisms made in Theorem 4.2, one of the main results of this paper,
can find its root in the classical Gr\"otzsch's theorem which is proved
by a length-area argument \cite{\Li}. Teichm\"uller generalized this
result to closed Riemann surfaces, in particular tori, by an ergodic
version of the length-area argument \cite{\Teich}. The notion of modulus
of a curve family is a formalism of the length-area (volume) argument.
In section 4, we reformulate Teichm\"uller's method to the CR setting
by computing moduli of some special families of curves and successfully
find the extremal quasiconformal homeomorphisms in certain homotopy
classes of mappings between CR circle bundles over flat tori.

\subhead Acknowledgement \endsubhead
This author is very grateful to his academic advisor, L\'aszl\'o Lempert
for the guidance with great insights. Thanks also to David Drasin and
Juha Heinonen for helpful talks.

\heading
2. Regularity of conformal homeomorphisms
\endheading

For $j= 1,2$, let $M_j$ be a smooth strongly pseudoconvex CR 3-manifolds.
The contact bundles $HM_j$ is assumed to be smooth and orientable, that is,
there exists smooth global 1-form $\eta_j$ on $M_j$, which is called a
contact form, so that $HM_j = \text{Ker }\eta_j$. Let $J_j: HM_j
\to HM_j$ denote the CR structure on $M_j$. $H^{0,1}M_j\triangleq\{X+iJ_jX\,
|\, X\in HM_j\}\i \Bbb C\otimes HM_j$ is the (0,1) tangent bundle on $M_j$.
$\wedge^{1,0}M_j \triangleq\{\text{ linear functional }\psi: \Bbb C\otimes
HM_j\to \Bbb C\, |\, \psi(J_j X) = i\psi (X),\text{ for } X\in HM_j\}$.

A mapping $f: M_1\to M_2$ is said to be absolutely continuous on lines (ACL)
if for any open set with a smooth contact fibration, $f$ is absolutely
continuous along all fibers in this fibration except a subfamily of measure
zero. Here a subfamily of fibers of the fibration is said to have measure zero
if intersections of these fibers with any transversal regular surface has
measure zero on the surface (see [12]).

\proclaim {Definition 2.1} A homeomorphism $f:M_1\to M_2$ is said to be
$K$-quasiconformal for a finite constant $K\ge 1$ if

(i) $f$ is ACL;

(ii) $f$ is differentiable almost everywhere, and its differential $f_*$
preserves the contact structures, i.e., $f_*(H_qM_1)\i H_{f(q)}M_2$, for
$q\in M_1$ where $f$ has differential.

(iii) for norms $|\cdot |_1$ and $|\cdot |_2$ defined by any Hermitian metrics
on $HM_1$ and $HM_2$ respectively,
$$
K(f)(q)=\frac{\underset{X\in H_qM_1, |X|_1=1}\to{\max} |f_*(X)|_2}
{\underset{X\in H_qM_1, |X|_1=1}\to{\min}|f_*(X)|_2} \le K<\infty ,
\tag2.1
$$
for almost all $q\in M_1$. Here $K=\underset{q\in M_1} \to{\text{ess sup}}
\, K(f)(q)$ is called the maximal dilatation of $f$. $f$ is conformal
if $K(f)=1$. $K(f)=\infty$ if $f$ is not $K$-quasiconformal for any finite
$K\ge 1$.
\endproclaim

A mapping $f: M_1\to M_2$ which is differentiable almost everywhere
is said to have $L^p_{\text{loc}}$ horizontal derivatives for $p\ge 1$ if
for any smooth function $h: M_2\to \Bbb R$ and any smooth local section
$X$ of $HM_1$ on an open set $U\i\i M_1$, the function $X(h\circ f)$ which
is defined almost everywhere on $U$ is in $L^p(U)$.

A mapping $f: M_1\to M_2$ is said to have $L^p_{\text{loc}}$ weak
horizontal derivatives for some $p\ge 1$ if for any smooth function $h: M_2\to
\Bbb R$ and any open set $U\i\i M_1$ with a smooth local section $X\ne 0$
of $HM_1$ on $U$, there exists a function $g\in L^p_{\text{loc}}(M_1)$
so that
$$
\int_U X\phi\cdot (h\circ f)\, dv_1 =-\int_U g\,\phi\, dv_1 \tag2.2
$$
for all $\phi\in C^\infty_0(U)$, where $dv_1$ is a smooth volume form
on $M_1$. Certainly, the function $g$ depends on the choice of $dv_1$.

For the proof of the next theorem, we fix a norm $|\,\cdot |$ on $HM_1$.
A regular curve on a contact manifold is called Legendrian if it is tangent
to the contact structure. For any smooth Legendrian curve $\gamma: I\to M_1$
with an interval $I\in \Bbb R$ and a function $g$ defined on an open
neighborhood of $\gamma$, define the line integral
$$
\int_\gamma g = \int_I g(\gamma(t))|\gamma^\prime (t)|dt.
$$

\proclaim {Theorem 2.2} A mapping $f: M_1\to M_2$ is ACL and has
$L^p_{\text{loc}}$ horizontal derivatives for some $p\ge 1$ if and only if
$f$ has $L^p_{\text{loc}}$ weak horizontal derivatives.
\endproclaim

\demo {Proof} First assume that the homeomorphism $f: M_1\to M_2$ is ACL
and has $L^p_{\text{loc}}$ horizontal derivatives. Let $h: M_2\to \Bbb R$
be a smooth function and $U\i\i M_1$ be any open set with a smooth section
$X\ne 0$ of $HM_1$ on it.

We can assume that the trajectories $\Gamma =\{\gamma\}$ of $X$ form a
contact fibration of $U$ by shrinking $U$ appropriately. Let $\Gamma_1$
be the subfamily of those $\gamma\in\Gamma$ along which $f$ is absolutely
continuous. Then $\Gamma\setminus \Gamma_1$ has measure zero. Along
$\gamma\in\Gamma_1$, $h\circ f$ is absolutely continuous, so $X(h\circ f)$
exists almost everywhere on $\gamma$ and
$$
\int_\gamma X\phi\cdot (h\circ f) +\int_\gamma \phi\cdot X(h\circ f) 
=\int_\gamma X(\phi\cdot h\circ f) =0,\qquad \forall \phi\in C_0^\infty(U).
\tag2.3
$$

We have topological and differential structures on $\Gamma$ such
that the natural projection $p: U\to\Gamma$ is open
and smooth. Then $\Gamma$ becomes a smooth surface. Let $t$ be the parameter
of the flow generated by $X$ and $\omega$ be any area form on $\Gamma$.
Then $dv_1 \triangleq p^*\omega\wedge dt$ is a volume form of $U$. Integrating
the expressions in (2.3) against $\omega$ with respect to $\gamma\in\Gamma_1$,
then by the ACL property of $f$, local $L^p$ integrability of $X(h\circ f)$ and
Fubini's theorem, we obtain
$$
\int_U X\phi\cdot (h\circ f)\,dv_1 = - \int_U \phi\cdot X(h\circ f)\,dv_1.
\tag2.4
$$
So the weak derivative of $h\circ f$ in the $X$ direction is given by
$X(h\circ f)\in L^p_{\text{loc}}$.

Conversely, assume $f$ has $L^p_{\text{loc}}$ weak horizontal derivatives.
For any open set $U\i\i M_1$ and a smooth contact fibration $\Gamma$ of
$U$, let $X$ be the nonzero horizontal vector field on $U$ so that $X_q$,
for any $q\in U$, is the tangent vector at $q$ of the fiber $\gamma\in \Gamma$
passing through $q$. Let $B\i\i U$ be an open set with the coordinate system
$\{(x,y,t)\,|\, a_1<x<a_2,\, b_1<y<b_2,\, c_1<t<c_2\}$, here $t$ is the
parameter of the flow generated by $X$, i.e., $X={\partial \over{\partial t}}$.

For any smooth function $h$ on $M_1$, $h\circ f$ has $L^p(U)$ weak derivative
in the direction $X$. Denote it by $\psi \in L^p(U)$. Hence there
exists a sequence of $C^1$ functions $g_n$ on $B$ such that $g_n$ converges
to $h\circ f$ uniformly in $B$ and $Xg_n$ converges to $\psi$ in
$L^p(B)$. Let $B_{x,y,t}=(a_1, x)\times (b_1,y)\times (c_1,t)$, $R_{x,y}
=(a_1,x)\times (b_1,y)$.
$$
\int_{B_{x,y,t}}Xg_n(u,v,w)\,dudvdw=\int_{R_{x,y}}(g_n(u,v,t)-g_n(u,v,c_1))
\,dudv.\tag2.5
$$
Hence by taking limits, we have
$$
\int_{B_{x,y,t}}\psi (u,v,w)\,dudvdw=\int_{R_{x,y}}((h\circ f)(u,v,t)
-(h\circ f)(u,v,c_1))\,dudv.\tag2.6
$$
Let $\{ t_n\}$ be a countable dense set of $(c_1, c_2)$. (2.6) implies
for each $t_n$, there exists a set $E_n\i R_{a_2,b_2}$ so that
$R_{a_2,b_2}\setminus E_n$ is of measure zero and
$$
\int^{t_n}_{c_1} \psi (x,y,w)dw=(h\circ f)(x,y,t_n)-(h\circ f)(x,y,c_1),
\qquad\forall\, (x,y)\in E_n. \tag2.7
$$
Then (2.7) is true for all $(x,y)\in E\triangleq\cap E_n$ and all $t_n$.
Hence by continuity of both sides in $t$,
$$
\int^t_{c_1} \psi (x,y,w)dw=(h\circ f)(x,y,t)-(h\circ f)(x,y,c_1),
\qquad\forall\, (x,y)\in E,\, t\in (c_1,c_2). \tag2.8
$$
So $(h\circ f)(x,y,t)$ is absolutely continuous in $t$ for $(x,y)\in E$.
Note $R_{a_2,b_2}\setminus E$ has measure zero. So $f$ is ACL since
$B\i\i U$ is an arbitrary rectangular coordinate chart and $U\i\i M_1$ is
arbitrary in $M_1$. Moreover, $X(h\circ f)$ exists almost everywhere
and  $X(h\circ f)=\psi\in L^p(U)$ on $U$. So $f$ has $L^p_{\text{loc}}$
horizontal derivatives. \qed
\enddemo

A homeomorphism $f: M_1\to M_2$ is said to be weakly CR if for any
smooth CR function $h: M_2\to \Bbb C$, open set $U\i\i M_1$, $\phi\in
C_0^\infty (U)$ and $\overline Z\in H^{0,1}M_1$, $\dsize\int_U\,
(h\circ f)\cdot\overline Z\phi =0$.

\proclaim {Theorem 2.3} Assume $M_1$ and $M_2$ are two smooth, strongly
pseudoconvex, locally embeddable CR 3-manifolds, $f: M_1\to M_2$ is a conformal
homeomorphism with $L^1_{loc}$ horizontal derivatives. Then $f$ is smooth
and CR.
\endproclaim

\demo {Proof} A simple linear algebra argument shows that at a point
$q\in M_1$ where $f$ is differentiable
$$
K(f)(q)= \frac{1+ |\mu(q)|}{1 - |\mu(q)|} \qquad\text{with}\qquad
|\mu(q)| = \Big|\frac{\langle f^*\psi_2, \overline Z\rangle }{\langle
f^*\psi_2, Z\rangle }\Big|(q), \tag2.9
$$
for any nonzero $\overline Z\in H^{0,1}M_1$ and nonzero $\psi_2\in\wedge^{1,0}
M_2$.

Then $K(f)=1$ implies that if $f$ is differentiable at $q\in M_1$,
$\langle f^*\psi_2, \overline Z\rangle = 0$ for any $\overline Z \in
H^{0,1}_qM_1$ and $\psi_2\in\wedge^{1,0}_{f(q)}M_2$. In other words,
$f_*$ preserves the CR structures at the points where it is defined.
Thus $\overline Z (h\circ f) = 0$ for any CR function $h$ on $M_2$ and
$\overline Z \in H^{0,1}M_2$. Theorem 2.2 says that $f$ is weakly CR.
$M_1$ and $M_2$ are locally embeddable implies they are locally embeddable
into ${\Bbb C}^2$ as hypersurfaces. A theorem of Pinchuk and Tsyganov
(Theorem 2, \cite {\Pinch}) asserts that such $f$ must be smooth, hence CR.
\qed
\enddemo

\heading
3. Moduli of curve families 
\endheading

Let $M$ be a smooth, compact, contact 3-manifold. We always assume $HM$
is smooth and oriented. A sub-Riemannian metric on $M$ with respect
to $HM$ is a smooth positive definite quadratic form on $HM$. Fix a
sub-Riemannian metric on $M$ with respect to $HM$ momentarily, and denote
by $|\cdot |$ the corresponding norm on $HM$. For general theory of
sub-Riemannian geometry, we refer to \cite{\Strich} and \cite{\Hamen}.

The sub-Riemannian metric on $M$ can be extended to a Riemannian metric
on $M$ canonically as follows. Let $\omega$ be the oriented area form
on $HM$ with respect to the sub-Riemannian metric. Then there exists
a unique contact form $\eta$ so that $d\eta |_{HM}=\omega$. Let $T$ be the
characteristic vector field of $\eta$, namely, T is the unique vector
field satisfying that $T\lrcorner\eta =1$ and $T\lrcorner d\eta =0$.
Declaring $T$ is a unit vector orthogonal to $HM$, we obtain a Riemannian
metric which is called the canonical extension of the sub-Riemannian metric.
The positive volume form of this Riemannian metric is $dv\triangleq
d\eta\wedge\eta$.

A curve $\gamma: I_\gamma \to M$ with an interval $I_\gamma\i\Bbb R$
is called locally rectifiable if $\gamma$ is absolutely continuous
and $\gamma^\prime (t)$ is tangent to $HM$ for almost all $t\in I_\gamma$.
$\gamma$ is called rectifiable if $\gamma$ is locally rectifiable and
the length
$$
l(\gamma )\triangleq\int_{I_\gamma}|\gamma^\prime (t)|\, dt <\infty. \tag3.1
$$
We set $l(\gamma )=\infty$ if $\gamma$ is not rectifiable.
For a locally rectifiable curve $\gamma$ and a non-negative
Borel-measurable function $\sigma$ on $M$, define the line integral
$$
\int_\gamma \sigma=\int_{I_\gamma}\sigma (\gamma (t))\,|\gamma^\prime
(t)|\,dt. \tag3.2
$$

\proclaim {Definition 3.1} Let $\Gamma$ be a family of curves $\gamma:
I_\gamma\to M$. An admissible measure for $\Gamma$ is a Borel-measurable
function $\sigma :M\to\Bbb R$ so that $\sigma \ge 0$ and $\int_\gamma
\sigma\ge 1$, for all locally rectifiable $\gamma\in\Gamma$. Denote
the set of admissible measures for $\Gamma$ by $A(\Gamma )$. The modulus
of $\Gamma$ is defined by
$$
\text{Mod}_M(\Gamma)=\underset{\sigma\in A(\Gamma )}\to{\inf}\int_M \sigma^4 dv.
\tag3.3
$$
\endproclaim

\demo {Remark} (1) It is easy to see that if two sub-Riemannian metrics on
$M$ with respect to $HM$ define the same conformal structure on $HM$, then
they give the same value to $\text{Mod}_M(\Gamma )$.

(2) If $\Gamma_r\i\Gamma$ consisting of all locally rectifiable curves
of $\Gamma$, then $\text{Mod}_M(\Gamma_r)=\text{Mod}_M(\Gamma)$.
\enddemo

The following proposition shows that modulus, regarded as a measure of
(more precisely, locally rectifiable) curve families, generalizes the
concept of measure zero used in the definition of ACL property. Thereafter
if a property holds for all curves in a family $\Gamma$ except a subfamily
with zero modulus, we say this property is true for almost all curves in
$\Gamma$. 

\proclaim {Proposition 3.2} Let $U$ be an open set of $M$, $\Gamma$ a
contact fibration of $U$, $\Gamma_1\i\Gamma$. Then $\Gamma_1$ has measure
zero if and only if $\text{Mod}_M(\Gamma_1)=0$.
\endproclaim

\demo {Proof} Without loss of generality, we assume $U$ is a domain of
the coordinator system $\{ (x,y,t)\,|\, a_1<x<a_2, b_1<y<b_2, c_1<t<c_2\}$
and $X={\partial\over{\partial t}}$ is tangent to $\Gamma$. Let $E\i
(a_1,a_2)\times (b_1,b_2)$ so that $\Gamma_1=\{\text{curves } t\mapsto
(x,y,t)\,|\, (x,y)\in E\}$.

If $\Gamma_1$ has measure zero, then $\int_Edxdy=0$. Notice 
$$
\sigma_0= \cases {1\over {\dsize c_2-c_1}},\qquad & \text{when } (x,y)\in E,
\, t\in (c_1,c_2), \\
0, \qquad & \text{otherwise},
\endcases
\qquad\in A(\Gamma_1). \tag3.4
$$
Therefore
$$
\text{Mod}_M(\Gamma_1)\le\, c\int_U\sigma_0^4\, dxdydt=0, \tag3.5
$$
where $c$ is a constant upper bound of the Jacobian $J$ on $U$ with
$J\,dxdydt=dv$.

If $\text{Mod}_M(\Gamma_1)=0$, then for any $\sigma\in A(\Gamma_1)$ with
$\sigma=0$ outside $U$, $\gamma\in\Gamma_1$,
$$
1\le (\int_\gamma\sigma )^4 \le (\int_\gamma\sigma^4)(\int_\gamma 1)^3
=(c_2-c_1)^3\,\int_\gamma\sigma^4. \tag3.6
$$
Taking the integral over $E$,
$$
\int_Edxdy\le \,c^\prime\,(c_2-c_1)^3\,\int_M\sigma^4. \tag3.7
$$
Hence
$$
\int_Edxdy\le \,c^\prime\,(c_2-c_1)^3\,\text{Mod}_M(\Gamma_1)=0, \tag3.8
$$
that is, $\Gamma_1$ has measure zero. \qed
\enddemo

Let $M_1$ and $M_2$ be two compact, smooth, strongly pseudoconvex CR
3-manifolds with smooth contact form $\eta_1$ and $\eta_2$ respectively.
Here the roles of sub-Riemannian metrics on $HM_1$ and $HM_2$ are played
by Hermitian metrics with respect to the CR structures $J_1$ and $J_2$
respectively.

\proclaim {Theorem 3.3} A $C^2$ homeomorphism $f: M_1 \to M_2$ is
$K$-quasiconformal for a constant $K>1$ if and only if for any family
$\Gamma$ of $C^1$ Legendrian curves on $M_1$
$$
{1\over K^2}\text{Mod}_{M_1}(\Gamma )\le \text{Mod}_{M_2}(f(\Gamma ))\le
K^2\text{Mod}_{M_1}(\Gamma ), \tag3.9
$$
where $f(\Gamma) =\{ f(\gamma )\,|\,\gamma\in\Gamma\}$.
\endproclaim

\demo{Proof} Assume $f: M_1\to M_2$ is $K$-quasiconformal. Then $f$ is contact, 
i.e., for contact forms $\eta_1$ and $\eta_2$ on $M_1$ and $M_2$ respectively,
$f^*\eta_2 = \lambda\,\eta_1$ with a $C^1$ function $\lambda$ on $M_1$. Then
$$
f^*(d\eta_2) = d\lambda\wedge\eta_1 + \lambda\, d\eta_1. \tag 3.10
$$
Therefore
$$
f^*(d\eta_2 |_{HM_2}) = \lambda\, d\eta_1|_{HM_1}. \tag3.11
$$
For $j=1,2$, the Levi form $L_j$ on $M_j$ is a bilinear form on $HM_j$
defined by
$$
L_j(X, Y) = d\eta_j (X, J_jY), \qquad \text{for }X, Y\in HM_j.
\tag3.12
$$
By replacing $\eta_j$ by $-\eta_j$, if necessary, we can always assume
that $L_j$ is positive definite. Hence $L_j$ is a Hermitian form on $HM_j$.
With respect to the Levi forms on $M_1$ and $M_2$, we define
$$
\alignedat1
\lambda_1(q) & =\underset{Y\in H_qM_1, |Y|_1=1}\to{\max}|f_*(Y)|_2, \\
\lambda_2(q) & =\underset{Y\in H_qM_1, |Y|_1=1}\to{\min}|f_*(Y)|_2.
\endalignedat
\tag3.13
$$
Then (3.11) implies that $|\lambda | = \lambda_1\lambda_2$. Moreover, (3.10)
implies that
$$
f^*(d\eta_2\wedge\eta_2) = \lambda^2 \,d\eta_1\wedge\eta_1. \tag3.14
$$
Thus the Jacobian of $f$ with respect to the volume forms $dv_1 =
d\eta_1\wedge\eta_1$ on $M_1$ and $dv_2 = d\eta_2\wedge\eta_2$ on $M_2$
is $J(f) = (\lambda_1\lambda_2)^2$.

For $\sigma_2\in A(f(\Gamma ))$,
$$
\alignedat1
\int_{f(\gamma )}\sigma_2 &  =\int_{I_\gamma}|f_*(\gamma^\prime(t))|_2\,
\sigma_2 (f(\gamma(t)))\,dt \\
& \le\int_{I_{\gamma}}\lambda_1\,|\gamma^\prime (t)|_1\,\sigma_2
(f(\gamma (t))) \,dt \\
& =\int_\gamma \lambda_1\cdot\sigma_2\circ f.
\endalignedat
\tag3.15
$$
Hence $\sigma_2\in A(f(\Gamma ))$ implies that $\lambda_1\cdot\sigma_2
\circ f\in A(\Gamma )$.

On the other hand,
$$
\alignedat1
\int_{f(\gamma )}\sigma_2 & \ge \int_{I_\gamma}\lambda_2\,|\gamma^\prime
(t)|_1\,\sigma_2(f(\gamma (t)))\,dt \\
& =\int_\gamma \lambda_2\cdot\sigma_2\circ f. 
\endalignedat
\tag3.16
$$
Hence $\sigma_1\in A(\Gamma )$ implies $(\sigma_1 /\lambda_2)\circ
f^{-1}\in A(f(\Gamma))$. Therefore,

$$
\alignedat1
\text{Mod}_{M_2}(f(\Gamma )) & =\underset{\sigma_2\in A(f(\Gamma))}\to{\inf}
\int_{M_2}\sigma_2^4 dv_2 \\
&  \le \underset{\sigma_1\in A(\Gamma )}\to{\inf}\int_{M_2}\big(
({\sigma_1\over\lambda_2})\circ f^{-1}\big)^4
dv_2 \\
& =\underset {\sigma_1\in A(\Gamma )}\to{\inf}\int_{M_1}{\sigma_1^4\over
\lambda_2^4}\cdot(\lambda_1\lambda_2)^2\,dv_1  \\
& \le K^2\,\text{Mod}_{M_1}(\Gamma ).
\endalignedat
\tag3.17
$$ 
$$
\alignedat1
\text{Mod}_{M_1}(\Gamma ) & =\underset{\sigma_1\in A(\Gamma )}\to{\inf}
\int_{M_1} \sigma_1^4\,dv_1 \\
& \le \underset{\sigma_2\in A(f(\Gamma ))}\to{\inf} \int_{M_1}(\lambda_1
\cdot\sigma_2 \circ f)^4\,dv_1 \\ 
& = \underset{\sigma_2\in A(f(\Gamma ))}\to{\inf}\int_{M_2}\frac{(\lambda_1
\circ f^{-1}\cdot\sigma_2 )^4}{(\lambda_1\lambda_2)^2\circ f^{-1}}\,dv_2 \\
& \le K^2\, \text{Mod}_{M_2}(f(\Gamma )).
\endalignedat
\tag3.18
$$

Assume that $f$ satisfies the inequalities (3.9). Then $f$ must preserve
the contact structures. Otherwise, there will be a point
$q\in M_1$ with a tangent vector $X_q\in H_qM_1$ so that
$f_*(X_q) \notin H_{f(q)}M_2$. $f$ is $C^1$, so there is a neighborhood
$U$ of $q$ with a nonzero smooth section $X$ of $HM_1$ on $U$ so that $f_*
(X_{q^\prime})$ is not contact for each $q^\prime\in U$. Let $\Gamma$
be the family of trajectories of $X$, and shrink $U$ appropriately such that
$\text{Mod}_{M_1}(\Gamma )\ne 0$. But each curve in $f(\Gamma )$ is not
Legendrian everywhere, hence not locally rectifiable. So $\text{Mod}_{M_2}
(f(\Gamma)) = 0$ by Definition 3.1. Then $f$ cannot satisfy (3.9) for this
family $\Gamma$. This contradiction shows that $f$ must be contact. 

If $f$ is not $K$-quasiconformal, there exists an open set $U\i M_1$
so that $\lambda_1 /\lambda_2 \ge K+\epsilon$, for some $\epsilon >0$.
Let $X$ be the vector field on $U$ so that $|X|_1=1,\, |f_*(X)|_2
=\lambda_2$. Let $\Gamma$ be the family of trajectories of $X$ in $U$.
Then for any integrable function $\sigma_2$ on $f(U)$,
$$
\int_{f(\gamma )}\sigma_2 = \int_\gamma \sigma_2\circ f\cdot |f_*(X)|_2
= \int_\gamma \lambda_2\cdot\sigma_2\circ f. \tag3.19
$$
So $\sigma_2\in A(f(\Gamma))$ if and only if $\sigma_1\triangleq\lambda_2
\cdot\sigma_2\circ f\in A(\Gamma)$. Hence
$$
\alignedat1
\text{Mod}_{M_2}(f(\Gamma)) & = \underset{\sigma_2\in A(f(\Gamma))}\to
{\inf}\int_{M_2}\sigma_2^4\, dv_2 \\
& = \underset{\sigma_1\in A(\Gamma)}\to {\inf}\int_{f(U)}\big(
({\sigma_1\over\lambda_2})\circ f^{-1}\big)^4\,dv_2 \\
& =\underset{\sigma_1\in A(\Gamma)}\to{\inf}\int_U ({\sigma_1 \over
\lambda_2})^4\cdot (\lambda_1\lambda_2)^2\, dv_1 \\
&  = \underset{\sigma_1\in A(\Gamma)} \to{\inf}\int_U ({\lambda_1\over
\lambda_2})^2\,\sigma_1^4 \,dv_1 \\
& \ge (K + \epsilon)^2\, \text{Mod}_{M_1}(\Gamma).
\endalignedat
\tag3.20
$$
Therefore, (3.9) implies that
$$
(K+\epsilon)^2\text{Mod}_{M_1}(\Gamma ) \le K^2 \text{Mod}_{M_1}(\Gamma ).
\tag3.21
$$
But we can always choose $U$ such that $\text{Mod}_{M_1}(\Gamma) \ne 0$.
So (3.9) cannot be true for such $\Gamma$. Hence $\lambda_1 /\lambda_2 \le K$
on $M_1$, i.e., $f$ is $K$-quasiconformal. \qed
\enddemo

We like to know if we can use the conclusion of Theorem 3.3 to define the
quasiconformality for a homeomorphism $f: M_1\to M_2$. The rest of this
section gives some partial solutions to this problem.

\proclaim {Theorem 3.4} If $f:M_1\to M_2$ is a homeomorphism so that
for a constant $K\ge 1$ and any curve family $\Gamma$ which forms a
smooth contact fibration of an open set in $M_1$
$$
\text{Mod}_{M_1}(\Gamma )\le K^2\,\text{Mod}_{M_2}(f(\Gamma )), \tag3.22
$$
then $f$ is ACL.
\endproclaim

\demo {Proof} Let $U\i M_1$ be an open set with a smooth contact fibration
$\Gamma$, $X\ne 0$ be a horizontal vector field on $U$ tangent to $\Gamma$.
By replacing $X$ by $X/ |X|$, we can assume $|X|=1$. Shrink $U$, if
necessary, so that there is a smooth surface $S\i U$ which intersects
each fiber of $\Gamma$ transversally once and only once. Parametrize
fibers $\gamma\in\Gamma$ by $t$ so that $X={\partial\over{\partial t}}$
and $\gamma (0)\in S$. $p:\,U\to S$ is the natural projection given by
$\gamma (t)\mapsto \gamma (0)$. 

Recall that Hermitian metrics on $HM_1$ and $HM_2$ can be extended canonically
to Riemannian metrics on $M_1$ and $M_2$ respectively. Restricting the
Riemannian metric on $M_1$ to $S$ makes $S$ a Riemannian 2-manifold.
Let $\omega$ be the area form on $S$, and $A_\omega (E)\triangleq \int_E\omega$
the $\omega$-area of a measurable set $E\i S$. Define a set function
$F$ for measurable set $E\i S$ by letting $F(E)=\text{vol}\,(f(p^{-1}(E)))$,
where vol refers to the Riemannian volume on $M_2$. Then Lebesgue's theorem
(Theorem 23.5, \cite{\Vais}) asserts that $F$ has finite derivatives at all
points in $S_1$ with respect to $\omega$-area, for a subset $S_1\i S$ with
$A_\omega (S\setminus S_1) =0$. We want to prove $f$ is absolutely continuous
along the fibers of $\Gamma$ passing through points in $S_1$.

For $q\in S_1$, let $D_r\i S\cap U$ be a disc centered at $q$ with radius $r$,
$\gamma_q: I\to U$ be the fiber of $\gamma$ passing through $q$. Take any
sequence $t_1^\prime, t_1^{\prime\prime}, t_2^\prime, t_2^{\prime\prime},...,
t_k^\prime, t_k^{\prime\prime}\in I$ such that 
$$
t_1^\prime <t_1^{\prime\prime} <t_2^\prime <t_2^{\prime\prime}
...<t_k^\prime <t_k^{\prime\prime}. \tag3.23
$$
Let $\Delta t_j \triangleq t_j^{\prime\prime}-t_j^\prime$. When $\underset
{1\le j\le k}\to{\max}\Delta t_j$ is small enough,
$$
d_{r,j}\le \Delta t_j \le 2\,d_{r,j}, \tag3.24
$$
where $d_{r,j}$ is the sub-Riemannian distance between the set $B_{r,j}^\prime
\triangleq \{\gamma (t_j^\prime)\,|\, \gamma\in\Gamma ,\, \gamma (0)\in D_r\}$
and the set $B_{r,j}^{\prime\prime}\triangleq \{\gamma (t_j^{\prime\prime})\,|\,
\gamma\in\Gamma ,\gamma (0)\in D_r \}$. This is true because fibers of
$\Gamma$ are Legendrian and the geodesic with respect to the sub-Riemannian
metric tangent to $\gamma$ at $\gamma (t^\prime_j)$ is locally a length
minimizing curve (Theorem 5.4, \cite{\Strich}). For $j=1,2,...,k$, denote
$$
R_{r,j} \triangleq \{\text{points }\gamma (t)\,|\,\gamma\in\Gamma,
\gamma (0)\in D_r, t_j^\prime\le t\le t_j^{\prime\prime}\}
$$
and let
$$
\Gamma_{r,j}\triangleq\{\text{curves }\gamma_j:t\mapsto\gamma (t)\,|\,
\gamma\in\Gamma, \gamma (0)\in D_r, t_j^\prime\le t\le t_j^{\prime\prime}\},
$$
a contact fibration of $R_{r,j}$. Then by (3.22),
$$
\text{Mod}_{M_1}(\Gamma_{r,j}) \le K^2\,\text{Mod}_{M_2}(f(\Gamma_{r,j})).
\tag3.25
$$

Next we use length-volume argument to give an estimate for $\text{Mod}_{M_1}
(\Gamma_{r,j})$. For $\sigma\in A(\Gamma_{r,j})$ and $\gamma_j\in
\Gamma_{r,j}$,
$$
1\le (\int_{\gamma_j}\sigma )^4 \le (\int_{\gamma_j}\sigma^4)
(\int_{\gamma_j}\,1)^3 =(\Delta t_j)^3\int_{\gamma_j}\sigma^4 . \tag3.26
$$
Integrating each term against $\omega$ over $D_r$, then
$$
A_\omega (D_r)\le c\, (\Delta t_j)^3\int_{R_{r,j}}\sigma^4 dv_1. \tag3.27
$$
where $c$ is constant, $dv_1$ is the volume form of the Riemannian metric
on $M_1$. By (3.3) and (3.24),
$$
A_\omega (D_r)\le 8c\,{d_{r,j}}^3\, \text{Mod}_{M_1}(\Gamma_{r,j}). \tag3.28
$$

On the other hand, let $\delta_{r,j}$ be the sub-Riemannian distance
between $f(B_{r,j}^\prime)$ and $f(B_{r,j}^{\prime\prime})$. Then ${1/
\delta_{r,j}}\in A(f(\Gamma_{r,j}))$. Thus
$$
\text{Mod}_{M_2}(f(\Gamma_{r,j}))\le {1\over\delta^4_{r,j}} \text{vol}
(f(R_{r,j})). \tag3.29
$$

Combining (3.25), (3.28) and (3.29), we have
$$
(A_\omega (D_r))^{1\over 4}\,\delta_{r,j}\le \, (8cK^2)^{1\over 4}\,
{d_{r,j}}^{3\over 4}\big(\text{vol}(f(R_{r,j}))\big)^{1\over 4}. \tag3.30
$$
Summing (3.30) over $j$,
$$
\alignedat1
& (A_\omega (D_r))^{1\over 4}\sum^k_{j=1}\delta_{r,j} \\
\le & \,(8cK^2)^{1\over 4}\,\sum^k_{j=1}{d_{r,j}}^{3\over 4}
\big(\text{vol}(f(R_{r,j}))\big)^{1\over 4} \\
\le & \,(8cK^2)^{1\over 4}\, (\sum^k_{j=1} d_{r,j})^{3\over 4}
\big(\sum^k_{j=1} \text{vol}(f(R_{r,j}))\big)^{1\over 4} \\
\le & \,(8cK^2)^{1\over 4}\, (\sum^k_{j=1} d_{r,j})^{3\over 4}
\big(\text{vol} (f(p^{-1}(D_r)))\big)^{1\over 4} .
\endalignedat
\tag 3.31
$$
So by (3.24),
$$
\sum^k_{j=1}\delta_{r,j}\le\,(8cK^2)^{1\over 4}\,(\sum^k_{j=1}
\Delta t_j)^{3\over 4} \big(\frac{F(D_r)}{A_\omega (D_r)}\big)^{1\over 4} .
\tag3.32
$$
Letting $r\to 0$,
$$
\sum^k_{j=1}\delta_j\le \,\big( 8cK^2\, F^\prime(q)\big)^{1\over 4}
(\sum^k_{j=1} \Delta t_j)^{3\over 4}, \tag3.33
$$
where $\delta_j\triangleq \underset{r\to 0}\to{\lim}\delta_{r,j}$, i.e.,
the sub-Riemannian distance between $f(\gamma_q(t_j))$ and $f(\gamma_q
(t_{j+1}))$. Hence $f$ is absolutely continuous along $\gamma_q$, for $q\in
S_1$. So $f$ is ACL. \qed
\enddemo

\proclaim {Corollary 3.5} If a homeomorphism $f:M_1\to M_2$ satisfies (3.22),
then $f$ has horizontal derivatives almost everywhere on $M_1$ in the sense
that for any open set $U\i M_1$ with a smooth section $X$ of $HM_1$ on $U$
and any smooth function on $f(U)$, $X(h\circ f)$ exists almost everywhere
on $U$.
\endproclaim

\heading
4. Extremal quasiconformal homeomorphisms
\endheading

In this section we will construct CR 3-manifolds $M_1$, $M_2$ and a
quasiconformal diffeomorphism $f_0: M_1\to M_2$ such that $K(f_0) \le
K(f)$, for any $C^2$ quasiconformal homeomorphism $f: M_1\to M_2$
which is homotopic to $f_0$. $M_1$ to be constructed is the quotient
of the 3-dimensional Heisenberg group over a lattice. So we start with
the Heisenberg group.  

The 3-dimensional Heisenberg group $\Bbb H^3$ is the space $\Bbb R^3$
endowed with the group structure defined by
$$
(x,y,t)\, (u,v,,s) = (x+u, y+v, t+s+2yu-2xv). \tag4.1
$$
The standard contact structure on $\Bbb H^3$ is given by the contact form
$$
\tilde\eta = -{1\over 2}ydx + {1\over 2}xdy + {1\over 4}dt. \tag4.2
$$
The contact bundle $H\Bbb H^3 = \text{Ker }\tilde\eta$ has two global sections
$$
\tilde X = {\partial\over{\partial x}} + 2y{\partial\over{\partial t}},
\qquad \tilde Y = {\partial\over{\partial y}} - 2x{\partial\over{\partial t}}
\tag 4.3
$$
which span $H\Bbb H^3$ everywhere and are invariant under the left group
translation. The standard CR structure on $\Bbb H^3$ is
$$
\tilde J: H\Bbb H^3\to H\Bbb H^3,\,\,\,  \tilde X\mapsto \tilde Y,\,
\tilde Y\mapsto -\tilde X. \tag4.4
$$
Note that $d\tilde\eta = dx\wedge dy$, thus the sub-Riemannian metric on
$H\Bbb H^3$ determined by the area form $d\tilde\eta$ and the complex structure
$\tilde J$, i.e., Levi form, is the one making $\{\tilde X, \tilde Y\}$
orthonormal. The canonical Riemannian extension of this metric is the
Euclidean metric on $\Bbb R^3$. 

Next we study the geodesics on $\Bbb H^3$ with respect to this sub-Riemannian
metric. Since the metric is invariant under the left group translation,
thus it suffices to study geodesics joining the origin $O= (0,0,0)$ and
a generic point $q$. Denote by $p$ the projection from $\Bbb H^3$ to the
horizontal plane $P\triangleq\{ t=0 \}$. It is easy to see that if $\gamma$ is a
rectifiable curve in $\Bbb H^3$ with respect to the sub-Riemannian metric,
$p(\gamma)$ is rectifiable with respect to the Euclidean metric on $\{ t=0\}$,
and the respective lengths of $\gamma$ and $p(\gamma )$ coincide. The
characterization of the geodesics given in the next theorem was first 
given by Kor\'anyi by studying Euler-Lagrange equations \cite{\Kgeo}. The
following proof, which is due to Lempert \cite{\Leme}, is a geometric one.

\proclaim {Proposition 4.1} On $\Bbb H^3$, a minimal geodesic to connect the
origin $O$ and $q\in \Bbb H^3$ is on either a straight line if $q\in P$,
or a helix whose projection to $P$ is a circle if $q\notin P$.
\endproclaim

\demo{Proof} The conclusion is obvious if $q\in P$. Next we assume
$q\notin P$. We first consider the case $p(q)= O$. i.e., $q=(0,0,t)$ for
some $t\ne 0$. Let
$$
\gamma : [0, l]\to \Bbb H^3, \qquad s\mapsto (x(s), y(s), t(s))  \tag4.5
$$
be any oriented rectifiable curve joining $O$ and $q$. Then
$$
-2yx^\prime +2xy^\prime +t^\prime =0, \tag4.6
$$
almost everywhere on $\gamma$ by (4.2). Thus
$$
t = \int_0^l t^\prime ds = \int_{p(\gamma )} 2ydx-2xdy = -4 \int_\Omega
dx\wedge dy, \tag4.7
$$
where $\Omega$ is the 2-chain on $P$ so that $\partial \Omega =
p(\gamma )$. Note length of $\gamma$ with respect to the sub-Riemannian
metric is equal to the length of $p(\gamma )$ on $P$ with
respect to Euclidean metric. Obviously, the length of $p(\gamma )$ is the
minimal only when $\Omega$ is a simply connected domain with fixed area
$|t|/4$. Furthermore, $p(\gamma )$ must be a circle if $\gamma$
is a length minimizing geodesic joining $O$ and $q$, by the isoperimetric
property on the Euclidean plane. 

If $p(q)\ne O$, the above reasoning with some slight modifications can be
applied. For instance, $\Omega$ is a 2-chain bounded by the union of
$p(\gamma )$ and the straight line segment from $p(q)$ to $O$. \qed

\enddemo

Denote $A=(1,0,0), B=(\sigma ,\tau ,0), C=(0,0,1)\in \Bbb H^3$ with $\tau\ne 0$.
Define $\Gamma = \{ n_1A + n_2B + n_3C \, |\, n_1, n_2, n_3 \in \Bbb Z \}$.
The Euclidean translations by elements of $\Gamma$ define a $\Gamma$-action.
The contact structure and CR structure on $\Bbb H^3$ are invariant under
this $\Gamma$-action. We consider the quotient space $M_1\triangleq \Bbb H^3 /
\Gamma$. $M_1$ has a contact structure with the contact form $\eta$
satisfying $\pi^*\eta =\tilde\eta$ and a CR structure $J_1$ inherented
from $\tilde J$. Hence $M_1$ becomes a CR manifold. Let $\pi : \Bbb H^3
\to M_1$ denote the natural projection and $X= \pi_*\tilde X, Y= \pi_*\tilde
Y$. Then with respect to the sub-Riemannian metric uniquely determined
by $\eta$ and $J_1$, $\{ X, Y\}$ is orthonormal. Note also $\overline
Z_1= {1\over 2} (X+iY)\in T^{0,1}M_1$. There is a canonical $\Bbb R$-action
on $\Bbb H^3$ defined by
$$
(x,y,t)\mapsto (x, y, t+t^\prime ), \qquad \text{for } t^\prime \in \Bbb R.
\tag4.8
$$
This $\Bbb R$-action is free and preserves both contact and CR structures on
$\Bbb H^3$. The plane $P$ can be regarded as the quotient space of this
$\Bbb R$-action. The natural projection $p: \Bbb H^3\to P$ is isometric
between the sub-Riemannian metric on $\Bbb H^3$ and the Euclidean metric
on $P$. The free $\Bbb R$-action on $\Bbb H^3$ commutes with the
$\Gamma$-action. Thus it induces a free $S^1$-actions on $M_1$ so that $M_1$
becomes a circle bundle over the torus $S_1\triangleq \pi (P)$. The torus $S_1$
has a flat Riemannian metric inherented from the Euclidean metric on $P$.
Denote the natural projection of this circle bundle by $\pi_1: M_1\to S_1$.
Note the $S^1$-action preserves the contact and CR structures on $M_1$.
We call such a bundle $M_1\overset{\pi_1}\to{\to} S_1$ a CR circle
bundle over a flat torus.

For any constant $K\ge 1$, define a new CR structure $M_2$ on the contact
structure of $M_1$ by declaring that
$$
\overline Z_2 \triangleq {1\over 2} (\sqrt{K}X+{i\over\sqrt{K}}Y)
= \frac {K+1}{4\sqrt K} (\overline Z_1 +\frac{K-1}{K+1}Z_1) \tag4.9
$$
is an (0,1) tangent vector. Obviously, the identity mapping $f_0: M_1\to M_2$
is quasiconformal and $K(f_0)(q)=K$ for all $q\in M_1$. 

\proclaim {Theorem 4.2} Let $f: M_1\to M_2$ be a $C^2$ homeomorphism homotopic
to $f_0: M_1\to M_2$. Then $K(f)\ge K$.
\endproclaim 

Note the CR structure on $M_2$ is also invariant under the circle action. So
if $S_2$ denotes the same smooth torus as $S_1$, but endowed with the complex
structure induced from CR structure on $M_2$, then $M_2\overset{\pi}\to{\to}
S_2$ is also a CR circle bundle.

\proclaim {Lemma 4.3 (Strichartz)} For any $q$ in a sub-Riemannian manifold
$M$, there is an $\epsilon >0$ so that if $q_1\in M$ with $d(q_1, q)\le
\epsilon$, there exists a length minimizing curve joining $q_1$ and $q$.
\endproclaim

This is Lemma 3.2 in \cite{\Strich}, a proof was given there. If $\alpha$
is a curve, denote by $[\alpha ]$ the homotopy class with fixed end points
of $\alpha$. Let $l(\alpha )$ be the length of $\alpha$ with respect to
the sub-Riemannian metric on $M$.

\proclaim {Lemma 4.4} Let  $\alpha: [0,1]\to M$ be a curve on a compact
sub-Riemannian manifold to connect two points $q_1$ and $q_2$, $\inf
l(\beta )$ be taken over all $\beta\in [\alpha ]$. Then this infimum
is attained at a rectifiable curve $\tilde\alpha\in [\alpha ]$.
\endproclaim 

\demo {Proof} Let $\epsilon_q$ be the maximal $\epsilon >0$ for $q\in M$
determined by Lemma 4.3. Since $M$ is compact, $\underset{q\in M}\to
{\inf}\epsilon_q >0$. Take any $r$ so that $0<3r < \underset
{q\in M}\to{\inf}\epsilon_q$. Cover $M$ with finitely many closed
$r$-balls $B(q_j, r),\, j= 1, 2, ..., k$, with respect to the sub-Riemannian
metric. Note any two points in the same $B(q_j, r)$ can be joined by a
length minimizing curve within $B(q_j, 3r)$.

Let $L={\inf}\,l(\beta )$ over all ${\beta\in [\alpha ]}$. Let $\beta_n
\in [\alpha ]$ be curves such that $\lim_{n\to\infty} l(\beta_n)=L$.
Divide each $\beta$ by ordered points $p_{n0}=q_1, p_{n1},
p_{n2}, ..., p_{nN}=q_2 \in\beta_n$ into ordered subcurves $\tau_{n1},
\tau_{n2}, ..., \tau_{nN}$ so that each $\tau_{nj}\i B(q_{m_{nj}}, r)$ for
some integer $m_{nj}$ with $1\le m_{nj}\le k$. Note this $N$ is uniform for
all $n$ by appropriately choosing $\beta_n$.

Let $\sigma_{nj}$ be a length minimizing curve in $B(q_{m_{nj}}, 3r)$
to join the end points of $\sigma_{nj}$ and $\sigma_n=\sigma_{n1}
\sigma_{n2}...\,\sigma_{nN}$. Then since $B(q_{m_{nj}}, 3r)$ is
simply connected, $\sigma_{nj}\in [\tau_{nj}]$, and hence $\sigma_n
\in [\beta_n]=[\alpha ]$. Furthermore
$$
L\le l(\sigma_n)\le\sum^N_{j=1}l(\sigma_{nj})\le\sum^N_{j=1}
l(\tau_{nj})=l(\beta_n). \tag4.10
$$
Therefore $\underset{n\to\infty}\to{\lim}l(\sigma_n)=L$.

Since $M$ is compact, there is a sequence $\{n_j\}\i\Bbb N$ so that for
each $j=1,2, ..., N-1$, $p_{nj}$ is convergent to some point $p_j$,
when $n_j\to\infty$. Note $p_{j-1}, p_j\in B(q_{m_j}, r)$ for some integer
$m_j$ with $1\le m_j\le k$. Connect $p_{j-1}, p_j$ by a length minimizing
curve $\alpha_j\i B(q_{m_j}, 3r)$ and let $\tilde\alpha=\alpha_1\alpha_2
...\alpha_N$. Then $\tilde\alpha\in [\alpha ]$ and
$$
l(\tilde\alpha )=\underset{n_i\to\infty}\to{\lim}l(\sigma_{n_i})=L. \qed
\tag4.11
$$
\enddemo

The next two lemmas are the key properties to establish the extremality
of $f_0$ described by Theorem 4.2.

\proclaim {Lemma 4.5} The flow of diffeomorphisms $h_t$ generated by $X$
on $M_1$ preserve the volume form $dv_1 = d\eta\wedge\eta$.
\endproclaim 

\demo {Proof} This is obvious since the the flow of diffeomorphisms $\tilde
h_s$ on $\Bbb H^3$ generated by $\tilde X$  are given by
$$
\tilde h_s : (x,y,t)\mapsto (x+s, y, t+2sy), \,\, s\in\Bbb R  \tag4.12
$$
which preserve the volume form $d\tilde \eta\wedge\tilde\eta =dx\wedge
dy\wedge dt$. \qed
\enddemo

We define a sub-Riemannian metric on $M_2$ so that $\{\dsize\sqrt{K}X,{Y
\over \sqrt{K}}\}$ is orthonormal. Note this sub-Riemannian metric induces the
same CR structure on $M_2$. Denote the corresponding curve length by $l_2$.

\proclaim {Lemma 4.6} Let $\gamma :\Bbb R\to M_2$ be an integral line
of the vector field $\sqrt{K}X$. Then for $[a, b]\i \Bbb R$ and with
respect to the sub-Riemannian metric on $M_2$, the curve
$\gamma |_{[a,b]}$ is the unique length minimizing curve in its
homotopy class of curves with end points $\gamma (a)$ and $\gamma (b)$.
\endproclaim

\demo{Proof} Lifting the sub-Riemannian metric on $M_2$ to its universal
covering $\Bbb H^3$, we obtain a new sub-Riemannian metric on $\Bbb H^3$.
Let us call this new sub-Riemannian manifold $\Bbb H^3_2$. Then on
$\Bbb H^3_2$, there is a characterization of geodesics similar to
Proposition 4.1. The trajectories of $\sqrt{K} \tilde X$ are exactly
those straight geodesics. Indeed, the sub-Riemannian metric is invariant
under the $\Bbb R$-action (4.8), hence a rectifiable curve $\gamma$
in $\Bbb H^3_2$ has the same length with $p(\gamma )$ on $P_2$. Here $P_2$
is the plane $P$ endowed with the flat Riemannian metric so that $\sqrt{K}
{\partial\over{\partial x}}, {1\over\sqrt{K}}{\partial\over{\partial y}}$
are orthonormal. Therefore the $x$-axis, which is the trajectory of $\sqrt{K}
\tilde X$ passing through $O$, is a geodesic without conjugate points on
it. Note $\sqrt{K}\tilde X$ and the sub-Riemannian metric are invariant
under the group left translation. Hence this lemma is true. \qed
\enddemo

\demo {Proof of Theorem 4.2} Let $f:M_1\to M_2$ be an $C^2$ quasiconformal
homeomorphism. For $q\in M_1$, let $\gamma_{q,a}:[-a, a]\to M_1$ be the
integral line of $X$ so that $\gamma_{q,a}(0)=q$. Let $\Gamma_a =
\{\gamma_{q,a}\,|\, q\in M_1 \}$. Then by Theorem 3.4, 
$$
\text{Mod}_{M_1}(\Gamma_a)\le K(f)^2\,\text{Mod}_{M_2}(f(\Gamma_a)). \tag4.13
$$
First note the length of $\gamma_{q,a}$ is $2a$, so $1/ 2a\in A(\Gamma_a)$ and
$$
\text{Mod}_{M_1}(\Gamma_a)\le\int_{M_1}({1\over {2a}})^4dv_1 =({1\over {2a}})^4
\text{vol}(M_1). \tag4.14
$$
For $\sigma\in A(\Gamma_a)$,
$$
1\le (\int_{\gamma_{q,a}}\sigma)^4 \le (\int_{\gamma_{q,a}} 1)^3\,
(\int_{\gamma_{q,a}}\sigma^4) =(2a)^3\int_{\gamma_{q,a}}\sigma^4. \tag4.15
$$
That implies
$$
\frac{1}{(2a)^3}\le \int_{\gamma_{q,a}}\sigma^4 =\int_{-a}^a (\sigma (
\gamma_{q,a}(s)))^4\,ds. \tag4.16
$$
Taking the integrals of both sides of (4.16) against $dv_1$ with respect to $q$
over $M_1$,
$$
\alignedat1
\frac{1}{(2a)^3}\text{vol}(M_1) & \le\int_{M_1}(\int_{-a}^a (\sigma (
\gamma_{q,a}(s)))^4\,ds\,) dv_1 \\
& = \int_{-a}^a (\int_{M_1}(\sigma (\gamma_{q,a}(s)))^4\,dv_1)\,ds \\
& = 2a\,\int_{M_1}\sigma^4\,dv_1.
\endalignedat
\tag4.17
$$
The last equality in (4.17) is due to Lemma 4.5, which implies that
$$
\int_{M_1}(\sigma (\gamma_{q,a}(s)))^4\,dv_1 =\int_{M_1}\sigma^4\, h_s^*(dv_1)
= \int_{M_1}\sigma^4\,dv_1 \tag4.18
$$
is independent of $s$. Therefore
$$
({1\over {2a}})^4\text{vol}(M_1)\le \int_{M_1}\sigma^4\,dv_1,\qquad \forall
\,\sigma\in A(\Gamma_a). \tag4.19
$$
Combining (4.14) and (4.19), we get
$$
\text{Mod}_{M_1}(\Gamma_a) =({1\over {2a}})^4\text{vol}(M_1). \tag4.20
$$ 

Next we estimate $\text{Mod}_{M_2}(f(\Gamma_a))$ in (4.13) for a quasiconformal
homeomorphism $f$ homotopic to $f_0: M_1\to M_2$. Since $f$ is homotopic
to $f_0$, there is a continuous map $H: [0,1]\times M_1\to M_2$ so that 
$$
H(0,q)=f_0(q), \qquad H(1,q)=f(q), \qquad \forall q\in M_1. \tag4.21
$$
Let $\alpha_q: [0,1]\to M_2$ be a curve given by $t\mapsto H(t,q)$.
Define $G(q)=\inf l_2(\beta )$ over all $\beta\in [\alpha_q]$. Then by
Lemma 4.4, $G(q)$ is attained at a rectifiable curve
$\tilde\alpha_q\in[\alpha_q]$.

Now we prove that $G(q)$ is continuous on $M_1$. When $q_1$ is close enough
to $q$ on $M_1$, $f_0(q_1)$ is close enough to $f_0(q)$ so that there exist
a length minimizing curve $\delta_1$ joining $f_0(q)$, $f_0(q_1)$ according
to Lemma 4.3. Let $\delta_2 = f(\delta_1 )$. Then
$$
\delta_1\tilde\alpha_{q_1}\delta_2^{-1}\in [\alpha_q], \qquad
\delta_1^{-1}\tilde\alpha_q\delta_2\in [\alpha_{q_1}]. \tag4.22
$$
This implies
$$
\alignedat1
& l_2(\tilde\alpha_q)  \le l_2(\delta_1) + l_2(\tilde\alpha_{q_1})
+ l_2(\delta_2^{-1}), \\
& l_2(\tilde\alpha_{q_1}) \le l_2(\delta_1^{-1}) + l_2(\tilde\alpha_q)
+ l_2(\delta_2).
\endalignedat
\tag4.23
$$ 
In other words,
$$
|G(q_1)-G(q)| \le l_2(\delta_1 ) + l_2(\delta_2 ) . \tag4.24
$$
Then there exists $A>0$ so that
$$
G(q) \le A <\infty, \qquad \forall\, q\in M_1, \tag4.25
$$
since $M_1$ is compact.

For a curve $\gamma: [a,b]\to M_1$, $\tilde\alpha_{\gamma (a)} f(\gamma )
\tilde\alpha^{-1}_{\gamma (b)}\in [f_0(\gamma )]$. Then by Lemma 4.6 and
(4.25),
$$
l_2(f_0(\gamma ))\le l_2 (\tilde\alpha_{\gamma (a)}) +
l_2(f(\gamma)) + l_2 (\tilde\alpha^{-1}_{\gamma (b)}) \le 2A +
l_2(f(\gamma )). \tag4.26
$$
Applying (4.26) to $\gamma =\gamma_{q,a} \in \Gamma_a$, we obtain
$$
2\sqrt{K}a\le 2A+l_2(f(\gamma_{q,a})). \tag4.27
$$
Hence $1/(2\sqrt{K}a-2A)\in A(f(\Gamma_a))$. So we have the following estimate.
$$
\text{Mod}_{M_2}(f(\Gamma_a)) \le\int_{M_2}(\frac{1}{2\sqrt{K}a-2A})^4\,dv_2
= (\frac{1}{2\sqrt{K}a-2A})^4 \text{vol}(M_2). \tag4.28
$$
Note $\text{vol}(M_1)=\text{vol}(M_2)\ne 0$, so by (4.13), (4.20) and (4.28),
we have
$$
(\sqrt{K}-\frac{A}{a})^2\le K(f). \tag4.29
$$
Letting $a\to\infty$, we get that $K\le K(f)$.  \qed
\enddemo

\demo{Remark} In horizontal directions, the extremal quasiconformal
homeomorphism $f_0$ behaves as a stretching by the constant factor
$\sqrt K$ along trajectories of $X$ and a compressing by the same
factor along trajectories of $JX$. The generator $T$ of the circle action
is transversal to the contact bundle. In this transversal direction $T$,
$f_0$ is equivariant under the circle action since it is simply the identity
mapping while $M_1$ and $M_2$ have the same circle action. 
\enddemo

\enddocument